# A maximum likelihood estimation of Lévy-driven stochastic systems for univariate and multivariate time series of observations

Babak M. S. Arani[1, *]

[1]Bernoulli Institute for Mathematics, Computer Science and Artificial Intelligence, University of Groningen, Groningen, The Netherlands.

**Abstract**. Literature is full of inference techniques developed to estimate the parameters of stochastic dynamical systems driven by the well-known Brownian noise. Such diffusion models are often inappropriate models to properly describe the dynamics reflected in many real-world data which are dominated by jump discontinuities of various sizes and frequencies. To account for the presence of jumps, jump-diffusion models are introduced and some inference techniques are developed. Jump-diffusion models are also inadequate models since they fail to reflect the frequent occurrence as well as the continuous spectrum of natural jumps. It is, therefore, crucial to depart from the classical stochastic systems like diffusion and jump-diffusion models and resort to stochastic systems where the regime of stochasticity is governed by the stochastic fluctuations of Lévy type. Reconstruction of Lévy-driven dynamical systems, however, has been a major challenge. The literature on the reconstruction of Lévy-driven systems is rather poor: there are few reconstruction algorithms developed which suffer from one or several problems such as being data-hungry, failing to provide a full reconstruction of noise parameters, tackling only some specific systems, failing to cope with multivariate data in practice, lacking proper validation mechanisms, and many more. This letter introduces a maximum likelihood estimation procedure which grants a full reconstruction of the system, requires less data, and its implementation for multivariate data is quite straightforward. To the best of our knowledge this contribution is the first to tackle all the mentioned shortcomings. We apply our algorithm to simulated data as well as an ice-core dataset spanning the last glaciation. In particular, we find new insights about the dynamics of the climate in the curse of the last glaciation which was not found in previous studies.

**Introduction**

In many applied fields of science, it is customary to employ dynamic models based on first principles or empirical findings to describe the behavior of a phenomenon. However, these phenomenological models typically offer only a qualitative understanding of underlying dynamics, relying on assigned parameter values derived from experience or qualitative observations. In contrast, a more rigorous scientific approach involves working in reverse: using observational data or measurements over time to estimate a stochastic system of differential equations. This process, known as *'system reconstruction'* has been the subject of study for approximately two and a half decades (Siegert *et al.* 1998; Siegert & Friedrich 2001). Put differently, system reconstruction aims to disentangle the deterministic forces from the stochastic ones whose collective interaction shapes the dynamics reflected in data.

Diffusion models, also referred to as Langevin models, represent one of the simplest stochastic models and have been a focus of active research in the past two and a half decades (Siegert *et al.* 1998; Friedrich *et al.* 2011b; Kleinhans 2012; Raischel *et al.* 2012; Tabar 2019). These models offer a 'minimal' description of underlying dynamical laws, and various parametric and non-parametric inference algorithms have been developed for their reconstruction. Diffusion models consist of two interacting components: a deterministic part representing internal forces in the absence of perturbations, and a stochastic component reflecting the collective effect of numerous unmeasured or unseen external forces operating at smaller scales and sizes (Friedrich et al., 2011b). The interplay between these opposing forces leads to self-organized patterns and forms due to Haken's synergistic approach (Haken, 2006). The following describes a diffusion model in technical terms.

---

[1] Babak M. S. Arani, m.shojaeiarani@gmail.com, m.shojaei.arani@rug.nl

$$dx = \mu(x)dt + \sigma(x)dW, \qquad (1)$$

Where $\mu(x)$ represents the deterministic processes as a function of state $x$ and is referred to as the 'drift' vector, while $\sigma(x)$ reflects the impact of stochastic fluctuations per state $x$ and is termed the 'diffusion' matrix. Here, $W$ is a matrix of independent Wiener processes, ensuring that the noise source, $dW$, is white (uncorrelated) and Gaussian. In ecological contexts, for instance, $\mu(x)$ may describe vegetation growth (as a function of biomass $x$), which can be perturbed by environmental stochasticity with a 'weight' per state $x$ represented by $\sigma(x)$. If the intensity of environmental perturbations is constant across the state space, i.e., $\sigma(x) = \sigma$ is constant, then the noise in (1) is termed 'additive'; otherwise, it is termed 'multiplicative'. Both the drift vector and diffusion matrix can, in many real-world problems, also depend on time t, i.e., $\mu = \mu(x,t)$ and $\sigma = \sigma(x,t)$, implying that the diffusion model in (1) is non-stationary. Non-stationarity, in simpler terms, refers to the variation of statistical properties of a stochastic system over time, and almost all real-world processes exhibit non-stationarity. However, since the reconstruction of non-stationary systems typically follows the same techniques and algorithms developed to tackle stationary systems (for instance, one can consider a moving window over shorter periods of time where stationarity can be assumed, reconstruct the system over each window, and finally interpolate the results), we therefore focus our attention on stationary systems and their reconstruction hereafter. The noise source in (1) is white, which is an idealization implying that the autocorrelation function of noise immediately drops, known as '$\delta$-correlated' noise. In realistic situations, however, the autocorrelation function of noise usually decays exponentially based on a finite correlation time, leading to what is termed 'colored' noise. To account for colored noise, it is necessary to embed the system into higher dimensions where a white noise can be assumed in the extended state space (Häunggi & Jung 1994). For instance, if (1) represents a one-dimensional system under colored noise evolving based on a simple Ornstein-Uhlenbeck mechanism, then the extended two-dimensional embedding would have white noise.

The rationale behind Gaussian noise in (1) is that noise reflects our inadequate knowledge over an infinite number of unmeasured degrees of freedom, whose aggregate effect, bearing in mind the Central Limit Theorem (CLT), can be assumed to be Gaussian. An important premise of the CLT is the assumption that all contributing degrees of freedom in noise have 'finite variance'. The tails of a Gaussian distribution decay exponentially, rendering the Brownian noise in a diffusion model (equation (1)) to be *'continuous'*. However, realistic processes often involve jump *'discontinuities'* of various sizes and rates. For example, in ecology, natural disasters and unexpected extreme events such as fires, floods, earthquakes, and large storms are omnipresent and constitute an integral part of environmental stochasticity. In such situations, the assumption of 'finite variance' for the unmeasured environmental variables in the CLT can be violated, questioning the use of Gaussian noise. To see the difference between a diffusion system in equation (1) and a stochastic system dominated by jumps, it is worth mentioning a consequence of an important theorem called the Pawula Theorem (Risken & Risken 1996). According to this theorem, in a diffusion model, the following Kramers-Moyal (KM) coefficients

$$D^{(n)}(x) = \lim_{\tau \to 0} K^{(n)}(x,\tau)/n!\,\tau, \quad n = 1,2. \qquad (2)$$

(Where $K^{(n)}(x,\tau) = \langle (x(t+\tau) - x(t))^n | x(t) = x \rangle$ is a conditional moment) of order higher than two vanish, i.e., $D^{(n)} \equiv 0$ for $n > 2$. In contrast, in a process characterized by jumps all KM coefficients are present. Indeed, KM coefficients of order higher than two contain information about the characteristics of jumps (Risken & Risken 1996; Tabar 2019). Therefore, to fully describe the nature of jumps in a system, it is necessary to consider all such high-order KM coefficients.

To address the presence of jump discontinuities jump-diffusion models have been proposed (Bandi & Nguyen 2003)

$$dx = \mu(x)dt + \sigma(x)dW + \xi\, dJ, \tag{3}$$

where $J(t)$ is a Poisson jump process with jump frequency $\lambda(x)$ and jump size $\xi$ which is typically assumed to follow a zero-mean Gaussian density with variance $\sigma_\xi^2(x)$. However, a drawback of jump-diffusion models is their limited ability to adequately account for the characteristics of jumps. To see this, note that the reconstruction of jump-diffusion models in (3) requires KM coefficients up to order 6 only:

$$\sigma_\xi^2(x,\tau) = \frac{M^{(6)}(x,\tau)}{5M^{(4)}(x,\tau)}, \lambda(x) = \frac{M^{(4)}(x,\tau)}{3\sigma_\xi^2(x,\tau)}, \mu(x) = \lim_{\tau\to 0} M^{(1)}(x,\tau)/\tau, \\ \sigma(x) = \lim_{\tau\to 0} M^{(2)}(x,\tau)/\tau - \lambda(x,\tau), \tag{4}$$

where $M^{(j)}(x), j = 1,\ldots,6$ in (4) are similar to the KM coefficients in (2) without the factorial in the denominator. It is worth noting that a more realistic methodology for modeling and reconstructing jump-diffusions integrates KM coefficients up to order 8 (see section 12.3 in (Tabar 2019) for details), although it is less commonly used. So, using jump-diffusion models we cannot account for the infinite number of higher order KM coefficients. Moreover, higher-order KM coefficients have higher uncertainty, and as the KM order increases, more data is needed to maintain the same level of accuracy, making the reconstruction of jump-diffusion systems highly data-hungry.

**Non-Poisson Lévy-driven stochastic systems**

As argued in the introduction, the justification for using well-known Brownian noise is rooted in the Central Limit Theorem (CLT), which assumes that all contributing drivers of noise have finite variance. However, this assumption is often violated for realistic disturbances, leading to the proposal of jump-diffusion models by the help of Poisson noise. Despite this advancement, we further argued that the use of Poisson noise is inadequate, as it fails to incorporate all the necessary information to describe the jumpy character of stochastic environmental fluctuations. Poisson jumps have two key features that make them inappropriate in describing natural humps: they are large and infrequent (Aıt-Sahalia 2004). In fact, a Poisson jump process is the only pure Lévy process where the jump frequency is 'finite' over a finite time interval (Aıt-Sahalia 2004). In the real world, however, jump sizes form a continuum ranging from small to medium to large, and jump rates occur frequently within a finite time interval. These limitations invalidate the use of jump-diffusion models as suitable models for describing natural jumps. Natural jumps are indeed rare events but not that rare the Brownian and Poissonian types of perturbations predict

Therefore, we need to adopt a modeling strategy that utilizes a generalized CLT to account for the continuous spectrum of noise amplitudes and the frequent occurrence of jump discontinuities over a finite interval of time. This leads us to consider the following non-Poisson Lévy-driven stochastic system (for now we focus on the one-dimensional case but later elaborate on the multidimensional case)

$$dx = \mu(x)dt + \sigma(x)dL. \tag{5}$$

In (5), $L$ represents a Lévy motion, and its increments follow an $\alpha$-stable distribution. Here, $\sigma(x)$ is termed the 'noise intensity' which differs from diffusion models where it was referred to as the diffusion function to emphasize that the noise in (5) is no longer diffusive. A stable distribution is characterized by four parameters: the index ($\alpha$), skewness ($\beta$), scale ($\gamma$), and shift ($\delta$). The index parameter $\alpha$ is the most important parameter and varies in the interval [0 2]. It reflects the jumpy nature of noise: generally, a smaller $\alpha$ indicates jumpier behavior in terms of both magnitude and frequency of jumps and vice versa. However, note that any jump amplitude is possible for any $\alpha < 2$, but the likelihood of larger jumps decreases as $\alpha$ increases and vice versa. At the extreme limit $\alpha = 2$, however, a stable distribution collapses to Gaussian and this is the only value of $\alpha$ where noise is diffusive and we recover the Brownian noise, i.e., $dL = dW$ in (5). Indeed, if $\alpha < 2$ the character of noise changes fundamentally from being diffusive (continuous) to being advective (discontinuous). The skewness parameter $\beta$ influences the direction of jumps and ranges from -1 to 1: a negative $\beta$ associates with a left-skewed noise density, a positive $\beta$ associates with a right-skewed noise density and, $\beta = 0$ corresponds with a symmetric Lévy noise. The scale parameter $\gamma$ serves as an analog to the standard deviation in Gaussian distributions, but it is important to note that a stable distribution does not have a finite variance except when $\alpha = 2$ (in which we recover the Gaussian). Therefore, in some literature, $\gamma$ is referred to as the dispersion parameter. Finally, the shift parameter $\delta$ plays the role of 'mean' for the stable distribution but $\delta$ is mean only if $\alpha > 1$, otherwise the distribution lacks a finite mean. However, in real-world applications it is often the case that $\alpha > 1$. For convenience, we consider a 'standard' Lévy noise in (5) resembling a standard Gaussian distribution. By standard Lévy noise we mean a Lévy variable with zero shift parameter ($\delta = 0$) and unit dispersion parameter ($\gamma = 1$). The influence of non-zero shift parameters can be incorporated into the drift term $\mu(x)$, while the effect of non-unit scale parameters can be accommodated by the noise intensity $\sigma(x)$ in (5). The Euler discretization of equation (5) is as follows:

$$x(i+1) = x(i) + \mu(x(i))\Delta + \Delta^{\frac{1}{\alpha}} \sigma(x(i))\eta, \qquad (6)$$

where $\eta$ follows a standard Lévy motion with $\alpha$ as its index parameter.

### A short overview about pros and cons of some existing inference techniques for diffusion, jump-diffusion and Lévy-driven stochastic models

Despite their simplicity, diffusion models are not always realistic representations of complex phenomena. However, they can still serve as initial guesses prior to performing more advanced inference techniques. The primary objective of reconstructing diffusion models is to estimate the functions $\mu(x)$ and $\sigma(x)$ in (1) from time series of observations, say, $\{x_0, x_1, ..., x_N\}$. There are several and different reconstruction procedures for diffusion models. The *'Langevin approach'* is a nonparametric reconstruction algorithm, easy to implement for univariate data and very popular (Siegert *et al.* 1998; Rinn *et al.* 2016). The Langevin approach is relatively simple and involves the following steps: first, the state space spanned by the range of data is divided into bins. Then, for each bin center, estimates for the drift vector and diffusion matrix are obtained by calculating the 'mean rate of change' and 'variance of rate of change', respectively. These results are then interpolated across the state space. Here, we give more technical details for the case of univariate data (the case of multivariate is similar but more involved). Depending on the data we have, the accessible state space (which is the data range) should be binned in a way that each bin receives enough data and that there are sufficient number of bins. Then the first ($n = 1$) and second ($n = 2$) KM coefficients in () will be estimated where in the conditional expectation $K^{(n)}(x, \tau)$, $x$ is a bin center and $\tau = \Delta, 2\Delta, 3\Delta, ...$ runs over the first few multiples of time step $\Delta$ (which is assumed to be small so a rather high-resolution data is needed). To estimate the limit in (2), for a fixed $x$, one usually solves a linear regression problem of $K^{(n)}(x, \tau)$ versus $\tau = \Delta, 2\Delta, 3\Delta, ...$ . The first KM coefficient is an estimate for the drift function $\mu(x)$ and the second KM coefficient would be an estimate for $\sqrt{2D^{(2)}(x)}$. Unfortunately, the Langevin approach is not easy to implement in higher dimensions and a reliable estimation demands quite a lot of data. Furthermore, the Langevin approach lacks an intrinsic validation mechanism. Validation techniques typically involve comparing the statistical properties of the estimated model with those of the data. This may include calculating one-time conditional probabilities, such as $p(x_1|x_0)$, for both the model and the data across several initial states $x_0$ within the state space, and comparing the outcomes using a distance measure. For more details see, for instance, (Siegert & Friedrich 2001; Rinn *et al.* 2016)).

Another inference technique for diffusion models is called *'Euler reconstruction'* which is a parametric inference technique based on parametric models $\mu(x;\theta)$ and $\sigma(x;\theta)$ ($\theta$ being the vector of parameters) for the drift vector and diffusion matrix, respectively. The underlying principle involves utilizing the short-time propagator or transition density $P(x(t+\Delta)|x(t))$, $\Delta$ being small, to formulate a maximum likelihood estimation (MLE) procedure where the objective function is the following sum of log-likelihoods

$$\ell(\theta) = \sum_{i=1}^{N} \ln\{p(x((i+1)\Delta)|x(i\Delta);\theta)\}, \qquad (7)$$

where $N$ is the number of data points, typically the first unconditional term $\ln(p(x_0))$ is ignored, as we did, in (7) as it has a negligible impact and $p \sim \mathcal{N}(x(i\Delta) + \mu(x(i\Delta);\theta)\Delta, \Delta\sigma^2(x(i\Delta);\theta))$. Unlike the Langevin approach Euler inference technique incorporates its own validation mechanism, i.e., the model with the highest sum of log-likelihoods is the fittest. On the other hand, a disadvantage of Euler reconstruction is that it demands solving an optimization problem while the Langevin approach is based on simple linear regression problems. Nevertheless, with careful selection of models $\mu(x;\theta)$ and $\sigma(x;\theta)$ it would be easy to solve the MLE. An elegant approach involves employing polynomial (or even rational) models for the drift and diffusion functions. This transforms the problem into a polynomial optimization problem, which can be solved through multivariate algebraic systems of equations. Fast solvers are available for finding all local maxima with good precision (see, for instance (Bates *et al.* 2016)), avoiding the pitfalls of typical optimization procedures that may become trapped in local maxima.

Both Langevin and Euler inference techniques require relatively high-resolution time series data. It is important to note that there is a time scale known as the 'relaxation time scale' $\tau_R$ (Friedrich *et al.* 2011a; Honisch & Friedrich 2011; Honisch *et al.* 2012) so that if the sampling time $\Delta$ is bigger than this time scale then all inference techniques fail due to problems termed 'finite time effects' in the literature (Friedrich *et al.* 2011a; Honisch *et al.* 2012). This arises because, in the limit of 'statistical independence' (i.e., $\Delta > \tau_R$) consecutive data points become independent, rendering the dynamics indiscernible in such extremely low-resolution data (Anteneodo & Queirós 2010). Therefore, prior to performing any inference algorithm it is important to verify whether the dataset's resolution falls within an acceptable regime, i.e., $\Delta \lesssim \tau_R$. The groundbreaking approach of Aït-Sahalia, paved the way to fill the gap for univariate (Aït-Sahalia 2002) and multivariate (Aït-Sahalia 2002) sparsely sampled data. When the sampling time $\Delta$ is approximately of the same order of magnitude as $\tau_R$ (i.e., $\Delta \sim \tau_R$) we should use the inference technique by Aït-Sahalia instead of Langevin or Euler reconstruction techniques. Aït-Sahalia's inference technique is parametric and the core idea is to construct an approximate convergent closed form expansion, using Hermite polynomials, for the transition density $p(x(t+\Delta)|x(t))$ (and hence, for the log-likelihood function) when $\Delta$ is not small. However, Aït-Sahalia's approach is computationally expensive, particularly when $\Delta$ is rather big ($\Delta \sim \tau_R$) and the mentioned expansion fails to converge when parameters deviate significantly from the optimal parameter values, resulting in a locally defined optimization problem.

With respect to jump-diffusion models there exists a non-parametric inference approach whose spirit is the same as the Langevin approach for diffusion models (see, for instance, the references (Bandi & Nguyen 2003; Anvari *et al.* 2016)). The approach first estimates the KM coefficients $M^{(j)}(x,\tau), j=1,\dots,6$ in (4).

Then, jump size $\sigma_\xi^2$ will be estimated using 6$^{th}$ and 4$^{th}$ KM coefficients. Next, jump rate $\lambda$ will be calculated using the 4$^{th}$ KM coefficient and the already estimated jump size $\sigma_\xi^2$. Then, the drift term $\mu(x)$ will be estimated via a linear regression of $M^{(1)}(x,\tau)$ values versus $\tau = \Delta, 2\Delta, 3\Delta, \ldots$ (for a fixed $x$ and first few multiples of $\Delta$). Finally, $\sigma(x)$ will be estimated using the jump rate $\lambda$ and a linear regression of $M^{(1)}(x,\tau)$ values versus $\tau = \Delta, 2\Delta, 3\Delta, \ldots$ (again for a fixed $x$ and first few multiples of $\Delta$). Unfortunately, this technique has the same shortcomings of the Langevin approach: it is a data-hungry method and its implementation to higher dimensions is not so easy. Alternatively, parametric inference techniques for jump-diffusion models based on MLE principle have been developed. The methodology outlined in (Li & Chen 2016) constitutes an MLE inference technique for jump-diffusions, akin to Aït-Sahalia's approach for diffusion models. It constructs an approximate log-likelihood via a convergent series expansion of transition probability, tailored for sparsely sampled data. Nonetheless, it encounters similar issues to Aït-Sahalia's approach: computational complexity escalates with larger sampling times $\Delta$, and the approximate log-likelihood function may fail to converge as the departure from the true, albeit unknown, parameter values increases.

Finally, a non-parametric inference technique for stochastic systems driven by Lévy noise, as described by (Siegert & Friedrich 2001), shares similarities with the KM-based technique elucidated for diffusion models. Like its counterpart, this methodology primarily suits univariate data (and inherits other problems we discussed for KM-based techniques). However, a notable limitation is that it only works when the skewness parameter $\beta$ is 0, i.e., the case of symmetric Lévy noise. Furthermore, it can only handle index parameters $\alpha > 1$, a condition often met by real-world data. Yet, its estimation of $\alpha$ proves unreliable, posing a significant challenge.

## An MLE inference technique for stochastic systems driven by Lévy noise

Here, we give a parametric MLE-based inference technique where the basic idea is from the reference (Friedrich *et al.* 2011a) but their final result is not correct. Here, we elaborate on the approach and provide a correct result. The approach involves determining the transition probability $P(x((i+1)\Delta)|x(i\Delta))$ for stochastic systems driven by Lévy noise when sampling time $\Delta$ is small. Consider the following multidimensional system

$$dX = M(x)dt + \Sigma(x)dL, \qquad (8)$$

Where $X$ is a $d$-dimensional state variable, $dL$ is a diagonal $d$-by-$d$ matrix of independent Lévy noise sources, $\mathbf{M}(\mathbf{X}) = (\mu_1(\mathbf{X}), \mu_2(\mathbf{X}), \ldots, \mu_d(\mathbf{X}))^T$ represents a $d$-dimensional vector of drift forces and $\mathbf{\Sigma}(\mathbf{X})$ represents a $d$-by-$d$ diagonal matrix containing the corresponding noise intensities $\sigma_i(\mathbf{X})$ ($i = 1, \ldots, d$) for the noise sources in the matrix $d\mathbf{L}$. Therefore, our approach only treats the case where $\mathbf{\Sigma}$ is diagonal. As such, our approach is specifically designed to handle cases where $\mathbf{\Sigma}$ is diagonal, which is common in many real-world applications. Consider the following Euler discretization of (8)

$$X(i+1) = X(i) + M(X(i))\Delta + D \odot \Sigma(X(i))\,\eta$$

Where $\boldsymbol{\eta}$ is a diagonal matrix of independent noise sources $\eta_i(\mathbf{X})$ ($i = 1, \ldots, d$), each having a univariate standard stable distribution with corresponding index and skewness parameters $(\alpha_i, \beta_i)$, $D$ is a diagonal matrix with diagonal elements $\Delta^{1/\alpha_i}, i = 1, \ldots, d$ (see the univariate discretization in (6)) and $\odot$ is a Hadamard (or, elementwise) matrix product. For convenience, we define $\mathbf{G} = \mathbf{D} \odot \mathbf{\Sigma}$ and consider the following Euler discretization

$$X(i+1) = X(i) + M(X(i))\Delta + G(X(i))\boldsymbol{\eta}.$$

Assuming the matrix $G$ is invertible (i.e., $\sigma_i > 0$ for all $i$ and this is often the case in applications) we get

$$\boldsymbol{\eta} = G^{-1}(X(i))[X(i+1) - X(i) - M(X(i))\Delta], \qquad (9)$$

Now, using (9) and Jacobian formula we recover the following short-time transition probability

$$p(X(i+1)|X(t)) = p(\boldsymbol{\eta})J(X(i)), \quad J = det\left(\frac{\partial \boldsymbol{\eta}}{\partial X(i+1)}\right)$$

where $J$ is the Jacobian and $p(\boldsymbol{\eta})$ is a multivariate standard stable distribution for the stable random variables $\eta_i$. By (9) we find that $J = det\left(G^{-1}(X(i))\right)$. This together with the fact that $\eta_1, \eta_2, \ldots, \eta_d$ are independent simplifies the previous relation to

$$p(X(i+1)|X(i)) = \frac{1}{det(G(X(i)))} p(\eta_1) p(\eta_2) \ldots p(\eta_d). \qquad (10)$$

Since $G = D \odot \Sigma$ equation (10) can be simplified to

$$p(X(i+1)|X(i)) = \Delta^{-(1/\alpha_1 + 1/\alpha_2 + \cdots + 1/\alpha_d)} p(\eta_1) p(\eta_2) \ldots p(\eta_d). \qquad (11)$$

Equation (11) is a short-time transition probability for (8) and is our central result. This equation provides the material we need in order to construct a MLE procedure as below

$$\boldsymbol{\ell}(\theta) = \sum_{i=1}^{N} \ln\{p(X((i+1)\Delta)|X(i\Delta); \theta)\}, \qquad (12)$$

where $\boldsymbol{\ell}$ is the log-likelihood function and $\theta$ is the vector of unknown parameters for a parametric model with drift vector $M(X; \theta)$ and matrix of noise intensities $\Sigma(X; \theta)$ in (8). Finally, in order to assess the uncertainty of the estimated parameters, say $\hat{\theta}$, we follow a standard approach by calculating the following Fisher information (FI) matrix

$$\mathcal{F} = -\mathrm{E}(\mathcal{H}(\theta)),$$

where $\mathcal{H}(\theta) = \partial^2 \boldsymbol{\ell}/\partial\theta\partial\theta^T$ is the Hessian matrix. The FI matrix can be approximated by the 'observed' Fisher information as $\mathcal{F} \sim -\mathcal{H}(\hat{\theta})$. Then variance-covariance matrix $\Sigma_1 = \mathcal{F}^{-1}$ can be estimated as $\Sigma_1 \sim [-\mathcal{H}(\hat{\theta})]^{-1}$. Eventually, the error for the $i^{th}$ parameter in terms of standard deviation, say $\sigma_{\theta_i}$ can be estimated by the square root of the $i^{th}$ diagonal element of $[-\mathcal{H}(\hat{\theta})]^{-1}$.

**Application of the approach to some simulated and real datasets**

*Example 1: Univariate simulated data*

We begin with a simulated univariate dataset. Figure 1 depicts a time series with $N = 10^5$ data points simulated from the following Lévy-driven Landau system using a time step $\Delta = 0.01$

$$dx = (ax - bx^3)dt + \sigma\, dL, \quad a = b = 1, \sigma = 0.3, \alpha = 1.6, \beta = 0.$$

Here, as proof of concept, we consider a parametric model with drift $\mu(x) = ax - bx^3$ and noise intensity $\sigma(x) = c$, so that our vector of parameters is $\theta = [a, b, c, \alpha, \beta]$. After solving the MLE problem in (12) we obtained the following estimates for the model parameters

$$\hat{a} \sim 0.8618, \quad \hat{b} \sim 0.8981, \quad \hat{c} \sim 0.3031, \quad \hat{\alpha} \sim 1.5918, \quad \hat{\beta} \sim 0.0005$$

With the following uncertainties (in terms of standard deviation)

$$\sigma_{\hat{a}} \sim 0.0161, \quad \sigma_{\hat{b}} \sim 0.0112, \quad \sigma_{\hat{c}} \sim 0.0025, \quad \sigma_{\hat{\alpha}} \sim 0.0049, \quad \sigma_{\hat{\beta}} \sim 0.0097.$$

It is interesting to note that when using this MLE technique, the uncertainty of noise parameters is normally higher than that of drift parameters. While we have not investigated a formal proof for this claim, it has been verified for diffusion models (see, for instance, (Sorensen 2007; Tang & Chen 2009; Chang & Chen 2011)). We speculate that a similar proof might apply to stochastic systems under Lévy noise.

In general, it is not straightforward how to choose a proper model before embarking on MLE. While there are effective modeling strategies for univariate data, such as spline modeling, extending these methodologies to higher dimensions presents additional complexities. With spline modeling, we employ cubic splines models for $\mu(x; \theta)$ and $\sigma(x; \theta)$ instead of 'typical' parametric models we already considered. Using this spline modeling methodology, we do not need to worry about how to choose a 'proper' model. Splines are flexible structures which enables us to restore the unknown nonlinearities in the functions $\mu(x)$ and $\sigma(x)$. Here, we consider a (rather sparse) mesh over the state space called 'knot sequence' and the model parameters correspond to the values of the $\mu(x)$ and $\sigma(x)$ functions at knots. While the data range typically serves as our state space, datasets with few data points near the borders of data range can adversely affect parameter estimation. This problem is particularly pronounced in systems with Lévy noise due to the presence of rare but big jumps. To tackle this issue, we sometimes consider a state space slightly smaller than data range for estimating $\mu(x)$ and $\sigma(x)$, while including all data points when estimating index and skewness parameters (see the red dashed lines in Figure 1).

Figure 1 illustrates a spline modeling approach applied to the univariate dataset in this example. Here, we have employed an equidistance knot sequence with 6 knots across the considered state space. The estimated model parameters (i.e., the estimated values of $\mu(x)$ and $\sigma(x)$ at knots) are

$$\widehat{\mu(x)} \text{ at knots:} \quad 5.9201 \quad 0.44234 \quad -0.31423 \quad 0.36424 \quad -1.1589 \quad -7.9797$$

$$\widehat{\sigma(x)} \text{ at knots:} \quad 0.30318 \quad 0.30797 \quad 0.3055 \quad 0.30863 \quad 0.31114 \quad 0.41594$$

$$\hat{\alpha} \sim 1.5855, \quad \hat{\beta} \sim 0.009050,$$

with the corresponding uncertainties

$$\sigma_{\widehat{\mu(x)}} \text{ at knots:} \quad 0.1795 \quad 0.0154 \quad 0.0230 \quad 0.0262 \quad 0.0373 \quad 0.4848$$

$$\sigma_{\widehat{\sigma(x)}} \text{ at knots:} \quad 0.0219 \quad 0.0023 \quad 0.0036 \quad 0.0040 \quad 0.0050 \quad 0.0580$$

$$\sigma_{\hat{\alpha}} \sim 0.0050, \quad \sigma_{\hat{\beta}} \sim \sim 0.0110$$

*Example 2: Multivariate simulated data*

In this example we infer the parameters of the following three species Lotka-Volterra competition model using a simulated dataset with $N = 10^5$ data points and a time step $\Delta = 0.01$.

$$dX_i = r_i X_i \left(1 - \sum_{j=1}^{3} a_{ij} X_j\right) dt + \sigma_i \, dL_i, \quad i = 1,2,3,$$

with the following model parameters

$$r = \begin{bmatrix} 1 \\ 2 \\ 3 \end{bmatrix}, \quad a = \begin{bmatrix} 0.06 & 0.02 & 0.04 \\ 0.02 & 0.08 & 0.02 \\ 0.02 & 0.04 & 0.1 \end{bmatrix}, \quad \sigma = \begin{bmatrix} 0.3 \\ 0.3 \\ 0.3 \end{bmatrix}, \quad \alpha = \begin{bmatrix} 1.7 \\ 1.8 \\ 1.9 \end{bmatrix}, \quad \beta = \begin{bmatrix} -0.1 \\ 0.1 \\ 0.3 \end{bmatrix}.$$

Under the considered model parameters all species coexist (in the deterministic settings) and reach their equilibrium level of $\bar{X} = a^{-1}[1\ 1\ 1]^T \sim [10.87\ 8.69\ 4.35]^T$. The estimated model parameters are

$$\hat{r} \sim \begin{bmatrix} 0.9938 \\ 2.0447 \\ 3.0172 \end{bmatrix}, \quad \hat{a} \sim \begin{bmatrix} 0.0604 & 0.0220 & 0.0350 \\ 0.0198 & 0.0813 & 0.0177 \\ 0.0211 & 0.0370 & 0.1030 \end{bmatrix}, \quad \hat{\sigma} \sim \begin{bmatrix} 0.2968 \\ 0.2988 \\ 0.2998 \end{bmatrix}, \hat{\alpha} \sim \begin{bmatrix} 1.7080 \\ 1.8035 \\ 1.8981 \end{bmatrix}, \hat{\beta} \sim \begin{bmatrix} -0.0903 \\ 0.1022 \\ 0.3216 \end{bmatrix},$$

with the following corresponding uncertainties

$$\sigma_{\hat{r}} \sim \begin{bmatrix} 0.0396 \\ 0.0543 \\ 0.1142 \end{bmatrix}, \quad \sigma_{\hat{a}} \sim \begin{bmatrix} 0.0015 & 0.0021 & 0.0020 \\ 0.0007 & 0.0011 & 0.0014 \\ 0.0010 & 0.0017 & 0.0043 \end{bmatrix}, \quad \sigma_{\hat{\sigma}} \sim \begin{bmatrix} 0.0020 \\ 0.0017 \\ 0.0013 \end{bmatrix}, \sigma_{\hat{\alpha}} \sim \begin{bmatrix} 0.0047 \\ 0.0043 \\ 0.0035 \end{bmatrix}, \sigma_{\hat{\beta}} \sim \begin{bmatrix} 0.0138 \\ 0.0181 \\ 0.0276 \end{bmatrix},$$

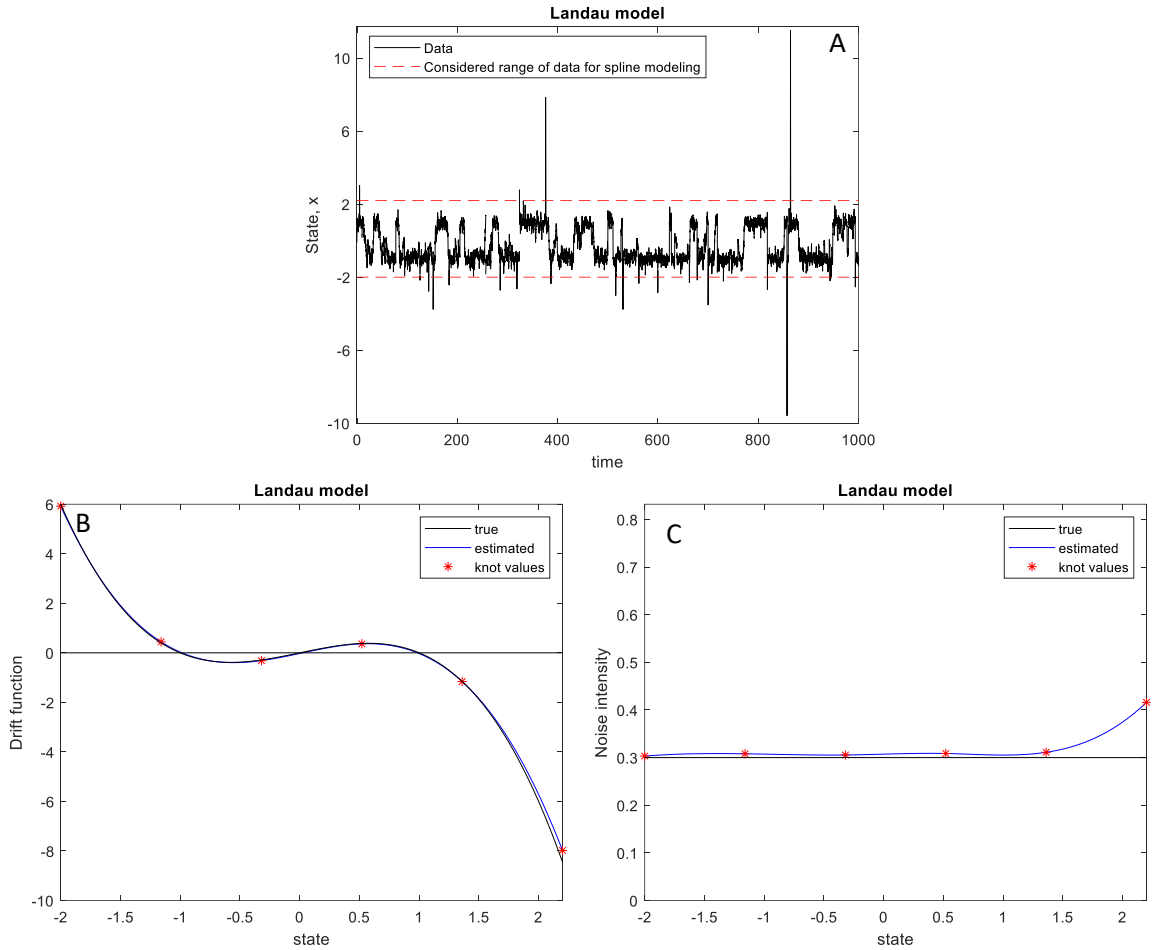

**Figure 1. Illustration of spline modeling for a univariate dataset simulated from Landau model.** In panel A, the black time series represents a dataset simulated by the Landau model under an additive Lévy noise using a time step of 0.01. The red dashed lines indicate the domain considered for $\mu(x)$ and $\sigma(x)$ in spline modeling. However, the entire dataset is taken into consideration when estimating index and skewness parameters. Panel B shows a spline model fitted to data in panel A. The black curve represents the true drift function, while the blue curve, a cubic spline, represents the estimated drift function. Red stars are 6 evenly spaced knots across the considered state space. In panel C the true noise intensity function is depicted by the black line, while the blue curve, a cubic spline, represents the estimated one. The red stars are the same as those in panel B.

## *Example 3: Univariate ice-core data*

Here, we reconstruct an ice-core calcium record from GRIP (Greenland Ice Core Project) (Fuhrer *et al.* 1993; Ditlevsen 1999). The resolution of this record which is nearly annual (Fuhrer *et al.* 1993). Additionally, it stands out as one of the longest records among glacial records, spanning from 11000 to 91000 years before the present (Fuhrer *et al.* 1993). This extensive time period covers most of the last glaciation, from the end of Eemian interglacial to the beginning of the current Holocene interglacial. Previous studies have demonstrated that the logarithm of this calcium record serves as a reliable climate proxy (Ditlevsen 1999) and is almost stationary (see Figure 2A).

This dataset has previously been reconstructed using various algorithms. For instance, the Langevin approach, discussed in detail in the third section, was employed in a study by (Arani *et al.* 2021). Other reconstruction procedures capable of handling stochastic systems driven by Lévy noise have also been employed (Ditlevsen 1999; Arani 2019). However, these algorithms are limited to symmetric Lévy noise, meaning they can only estimate the index parameter assuming the skewness parameter to be 0. Moreover, they may not reliably estimate the index parameter, as indicated by previous findings suggesting $\alpha \sim 1.75$. Nonetheless, this climate record exhibits significant jumps, which do not appear to be symmetric. Therefore, we have undertaken the new reconstruction in this paper to gain deeper insights into the nature of stochastic fluctuations during the curse of the last glaciation.

To infer the underlying data-generating system we fit a multiplicative stochastic system driven by Lévy noise to the logarithm of the calcium record. Utilizing spline models for $\mu(x)$ and $\sigma(x)$, as discussed in the first example, Figures 2B&C depict estimates of $\mu(x)$ and $\sigma(x)$, respectively, using 8 regularly spaced knots across the range [-2 3.5] (outside this range there is almost no data point). We find the estimate $\alpha \sim 1.58$ for the index parameter and the estimate $\beta \sim -0.2026$ for the skewness parameter. These estimates suggest that climate jumps during the last glaciation were more intense and frequent than previously thought. Additionally, these jumps were negatively skewed towards stadial states (see blue texts in Figures 2A&D). Figures 2B&C alone, do not provide clear evidence of alternative stable states in this climate record, a fact previously revealed by other studies (Babak M. S. Arani ; Ditlevsen 1999; Arani *et al.* 2021). To investigate this further, we calculate the '*effective potential*' given by $U(x) = -\log p(x)$ as shown in Figure 2D. This result aligns with the existing view that during the last glaciation the climate alternated between cold glacial regimes and warmer interstadial periods (Figure 2A&D), a phenomenon known as Dansgaard-Oeschger events (Dansgaard *et al.* 1993).

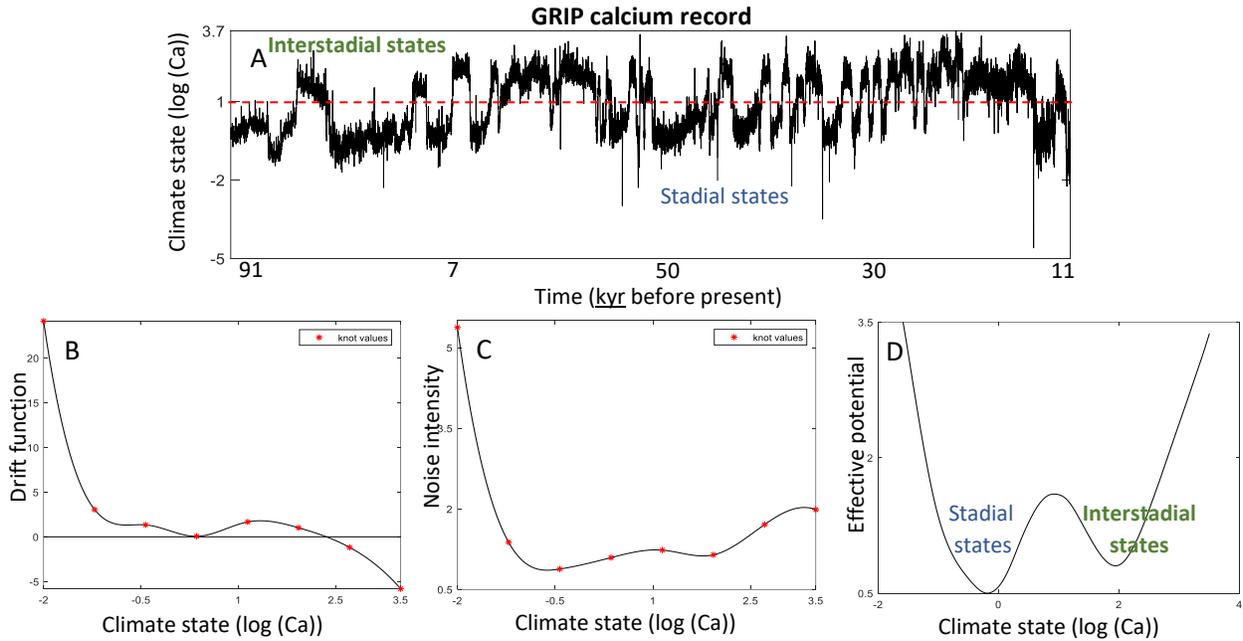

**Figure 2. Application of our inference technique to an ice-core climate record from the last glaciation.** The logarithm of calcium concentrations is shown to serve as a reliable climate proxy (A). The drift function $\mu(x)$ (B). Noise intensity $\sigma(x)$ (C). Effective potential $U(x) = -\log p(x)$ reveals the location of alternative climate regimes of cold interstadial and warmer interglacial states (D).

**Conclusion and future directions**

In this paper, we introduced a maximum likelihood inference technique to estimate the parameters of stochastic systems where the regime of stochasticity is of Lévy type. Previous algorithms suffered from several shortcomings, including being data-hungry, applicable only to specific classes of systems, and partially estimating noise parameters. Additionally, they often failed to tackle the multivariate case. However, our approach is only applicable when high-resolution data is available, presenting a major drawback shared by other existing inference algorithms. Addressing this gap for sparsely sampled data remains an open research area in the inference of Lévy-driven stochastic systems. Another challenge is that our methodology requires solving an optimization problem, which can become increasingly difficult as the number of parameters increases. This is due to the presence of a multitude of local maxima in such cases. While for diffusion systems of polynomial (or even rational) type it may be plausible to avoid solving the MLE using conventional optimization routines (which may get trapped in local maxima) by solving a multivariate algebraic system of equations instead, as the transition probability follows a Gaussian distribution, this approach may not be applicable to stochastic systems under Lévy noise. In these cases, the transition probability is not of Gaussian type. Nonetheless, further investigation into this idea is warranted. One possible avenue is to find an approximate series expansions of transition probability (see, for instance the ideas in (Ament & O'Neil 2018)) that could lead to a 'tamed' optimization problem of polynomial or rational type. This could potentially simplify the optimization process and improve computational efficiency.

# References

1. Aıt-Sahalia, Y. (2004). Disentangling diffusion from jumps. *Journal of Financial Economics*, 74, 487-528.

2. Aït-Sahalia, Y. (2002). Closed-form likelihood expansions for multivariate diffusions. National Bureau of Economic Research Cambridge, Mass., USA.

3. Aït-Sahalia, Y.J.E. (2002). Maximum likelihood estimation of discretely sampled diffusions: a closed-form approximation approach. 70, 223-262.

4. Ament, S. & O'Neil, M. (2018). Accurate and efficient numerical calculation of stable densities via optimized quadrature and asymptotics. *Statistics and Computing*, 28, 171-185.

5. Anteneodo, C. & Queirós, S.D. (2010). Low-sampling-rate Kramers-Moyal coefficients. *Physical Review E*, 82, 041122.

6. Anvari, M., Tabar, M.R.R., Peinke, J. & Lehnertz, K. (2016). Disentangling the stochastic behavior of complex time series. *Scientific reports*, 6, 35435.

7. Arani, B.M., Carpenter, S.R., Lahti, L., Van Nes, E.H. & Scheffer, M. (2021). Exit time as a measure of ecological resilience. *Science*, 372, eaay4895.

8. Arani, B.M.S. (2019). Inferring ecosystem states and quantifying their resilience: linking theories to ecological data. Wageningen University p. 112.

9. Babak M. S. Arani, S.R.C., Egbert H. van Nes, Ingrid A. van de Leemput, Chi Xu, Pedro G. Lind, and Marten Scheffer Stochastic regimes can hide the attractors in data, reconstruction algorithms can reveal them *Ecology*, Under revision.

10. Bandi, F.M. & Nguyen, T.H. (2003). On the functional estimation of jump–diffusion models. *Journal of Econometrics*, 116, 293-328.

11. Bates, D.J., Newell, A.J. & Niemerg, M. (2016). BertiniLab: A MATLAB interface for solving systems of polynomial equations. *Numerical Algorithms*, 71, 229-244.

12. Chang, J. & Chen, S.X. (2011). On the approximate maximum likelihood estimation for diffusion processes.

13. Dansgaard, W., Johnsen, S.J., Clausen, H.B., Dahl-Jensen, D., Gundestrup, N.S., Hammer, C.U. *et al.* (1993). Evidence for general instability of past climate from a 250-kyr ice-core record. *nature*, 364, 218-220.

14. Ditlevsen, P.D. (1999). Observation of α-stable noise induced millennial climate changes from an ice-core record. *Geophysical Research Letters*, 26, 1441-1444.

15. Friedrich, R., Peinke, J., Sahimi, M. & Tabar, M.R.R. (2011a). Approaching complexity by stochastic methods: From biological systems to turbulence. *Physics Reports*, 506, 87-162.

16. Friedrich, R., Peinke, J., Sahimi, M. & Tabar, M.R.R.J.P.R. (2011b). Approaching complexity by stochastic methods: From biological systems to turbulence. 506, 87-162.

17. Fuhrer, K., Neftel, A., Anklin, M. & Maggi, V. (1993). Continuous measurements of hydrogen peroxide, formaldehyde, calcium and ammonium concentrations along the new GRIP ice core from Summit, Central Greenland. *Atmospheric Environment. Part A. General Topics*, 27, 1873-1880.

18. Häunggi, P. & Jung, P. (1994). Colored noise in dynamical systems. *Advances in chemical physics*, 89, 239-326.

19. Honisch, C., Friedrich, R., Hörner, F. & Denz, C. (2012). Extended Kramers-Moyal analysis applied to optical trapping. *Physical Review E*, 86, 026702.

20. Honisch, C. & Friedrich, R.J.P.R.E. (2011). Estimation of Kramers-Moyal coefficients at low sampling rates. 83, 066701.

21. Kleinhans, D. (2012). Estimation of drift and diffusion functions from time series data: A maximum likelihood framework. *Physical Review E*, 85, 026705.

22. Li, C. & Chen, D. (2016). Estimating jump–diffusions using closed-form likelihood expansions. *Journal of Econometrics*, 195, 51-70.

23. Raischel, F., Russo, A., Haase, M., Kleinhans, D. & Lind, P.G. (2012). Searching for optimal variables in real multivariate stochastic data. *Physics Letters A*, 376, 2081-2089.

24. Rinn, P., Lind, P.G., Wächter, M. & Peinke, J.J.a.p.a. (2016). The Langevin Approach: An R Package for Modeling Markov Processes.

25. Risken, H. & Risken, H. (1996). *Fokker-planck equation*. Springer.

26. Siegert, S., Friedrich, R. & Peinke, J.J.a.p.c.-m. (1998). Analysis of data sets of stochastic systems.

27. Siegert, S. & Friedrich, R.J.P.R.E. (2001). Modeling of nonlinear Lévy processes by data analysis. 64, 041107.

28. Sorensen, M. (2007). Efficient estimation for ergodic diffusions sampled at high frequency. *CREATES Research Paper*.

29. Tabar, R. (2019). *Analysis and data-based reconstruction of complex nonlinear dynamical systems*. Springer.

30. Tang, C.Y. & Chen, S.X. (2009). Parameter estimation and bias correction for diffusion processes. *Journal of Econometrics*, 149, 65-81.